\documentclass[10pt,a4paper]{article}

\usepackage{graphicx}
\usepackage{amsmath,amsthm}
\usepackage{amssymb}
\usepackage[numbers]{natbib}
\usepackage{microtype}
\usepackage[final]{showkeys}

\numberwithin{equation}{section}
\allowdisplaybreaks

\newtheorem{Theorem}{Theorem}[section]
\newtheorem{Remark}{Remark}[section]
\newtheorem{Corollary}{Corollary}[section]

\newtheorem{Ex}{Example}[section]

\title{\textbf{Analytic solutions of fractional differential equations by operational methods}}

\author{$\text{Roberto Garra}_1$ \& $\text{Federico Polito}_2$\\
	\footnotesize (1) -- Dipartimento di Scienze di Base e Applicate per l'Ingegneria, ``Sapienza'' University, Rome\\
	\footnotesize Via A.\ Scarpa 16, 00161, Rome, Italy.\\
	\footnotesize Email address: rolinipame@yahoo.it \\
	\footnotesize (2) -- Dipartimento di Matematica, University of Torino\\
	\footnotesize Via Carlo Alberto 10, 10123, Torino, Italy.\\
	\footnotesize Email address: federico.polito@unito.it
	}

\begin{document}

	\maketitle

	\begin{abstract}

		\noindent We describe a general operational method that can be used in the analysis of fractional
		initial and boundary value problems with additional analytic conditions.
		As an example, we derive analytic solutions of some fractional generalisation of differential equations of
		mathematical physics. Fractionality is obtained by
		substituting the ordinary integer-order derivative with the Caputo fractional derivative.
		Furthermore, operational relations between ordinary and fractional
		differentiation are shown and discussed in detail.
		Finally, a last example concerns the application of the method to the study of
		a fractional Poisson process.
			
		\smallskip
	

	\end{abstract}
	
	\section{Introduction} 

		Recently, great interest has been devoted to the application
		of fractional calculus modelling to different fields of
		science e.g.\ rheology, biology, geomorphology (see for a recent review \citep{Rew,Podlb}).
		Indeed, it is well known that the introduction of fractional derivatives furnishes the system with
		a memory mechanism that is of great importance, in particular for some diffusive processes \citep{Benson}. 
		It is therefore important to understand the advantages of using fractional derivatives in classical equations.
		Currently, several different methods are used
		to solve fractional differential equations, from the classical Laplace transform
		for linear fractional equations \citep{Podlb}, to Adomian decomposition \citep{Adomian} and recently the
		Homotopy Perturbation Method \citep{Homotopy} which gives approximate solutions also for nonlinear
		fractional differential equations.
		In this note we provide analytic solutions of a general class of linear fractional differential equations by
		operational methods. In some previous papers, \citep{Dat1,Dat2} the authors investigated the differential
		isomorphism between Laguerre and ordinary derivatives. This approach permits us to study generalised
		evolution problems (see \citep{Dat4}). In our framework we generalise some of these
		results to fractional differential equations. 
		Besides that, we find a way to analytically solve boundary and initial value problems
		(BVP and IVP, respectively) of some interesting equations such as the space and time 
		fractional diffusion. Moreover, we present some applications considering widely used initial
		conditions to highlight the reliability	of our method.
		The organization of this paper is the following: in Section \ref{math} we discuss the relation between ordinary
		derivative and Caputo fractional derivative and
		discuss a stochastic interpretation of fractional differentiation from an operational point of view;
		in Section \ref{solution} we state the main result regarding the solution of 
		linear fractional differential equations by operational method and give some interesting examples. 

	\section{Relations between ordinary derivative and Caputo fractional derivative}
		
		\label{math}
		In \citep{Dat1}, the authors derived formal operational solutions of a class of
		partial differential equations involving the so-called Laguerre derivative operator
		\begin{align}
			\mathcal{D}_t = -\frac{\partial}{\partial t} t \frac{\partial}{\partial t},
		\end{align}
		in terms of Tricomi functions. In particular, for the initial value problem
		\begin{align}
			\begin{cases}
				\mathcal{D}_t f(x,t) = \mathcal{O}_x f(x,t), \\
				f(x,0) = g(x),
			\end{cases}
		\end{align}
		where $\mathcal{O}_x$ is a linear operator with
		constant coefficients on $x$ and $g(x)$ is an analytic function, the operational
		solution can be written as
		\begin{align}
			f(x,t) = C_0(t \mathcal{O}_x) g(x),
		\end{align}
		where
		\begin{align}
			C_0(y) = \sum_{n=0}^\infty \frac{(-y)^n}{(n!)^2}
		\end{align}
		is the zeroth order Tricomi function. Furthermore they pointed out the correspondence
		\begin{align}
			\mathcal{D}_t \longrightarrow \frac{\partial}{\partial t}, \qquad \qquad C_0(y) \longrightarrow
			e^y,
		\end{align}
		and also that if $h(x,t)$ is the classical solution, then $h(x, \Delta_t^{-1})$
		is the solution to the same problem but based on Laguerre derivative, where
		\begin{align}
			\Delta_t^{-1} \xi(t) = \int_0^t \xi(y) \mathrm dy.
		\end{align}
		
		In this paper, we exploit the same framework in order to arrive at solutions to fractional
		differential equations.
		The Caputo fractional derivative is defined as
		\begin{align}
			\label{caputo}
			D^{\nu} f(x) = \frac{1}{\Gamma(m-\nu)} \int_0^x \frac{\frac{\mathrm d^m}{
			\mathrm d y^m}f(y)}{(x-y)^{\nu-m+1}}
			\mathrm d y, \qquad \nu \in (0,+\infty), \: m=\lceil \nu \rceil,
		\end{align}
		and reduces to the classical derivative when $\nu=m$.
		It is well-known that the eigenfunction associated with the Caputo derivative is the Mittag--Leffler
		function
		\begin{align}
			E_\nu (x) = \sum_{r=0}^\infty \frac{x^r}{\Gamma(\nu r +1)}, \qquad \nu \in (0,+\infty), \:
			x \in \mathbb{R}.
		\end{align}
		This can be ascertained by simply writing for example the IVP problem
		\begin{align}
			\label{pagnotta}
			\begin{cases}
				D^\nu_t f(t) = \alpha f(t), & \nu \in (0,2], \\
				f(0) = 1, \\
				f'(0) = 0,
			\end{cases}
		\end{align}
		and by applying the Laplace transform to both members of \eqref{pagnotta}, thus obtaining
		\begin{align}
			z^\nu \mathcal{L} \{ f \} (z) - z^{\nu-1} = \alpha \mathcal{L} \{ f \} (z).
		\end{align}
		Now, by considering that $\mathcal{L} \{ f \}(z) = z^{\nu-1}/(z^\nu-\alpha)$ is the Laplace
		transform of a Mittag--Leffler function we readily arrive at
		\begin{align}
			f(t) = E_\nu(\alpha t^\nu), \qquad t \ge 0, \: \nu \in (0,2].
		\end{align}
		
		In order to complete the correspondence between classical and Caputo derivative, first
		we notice that the natural choice of the Riemann--Liouville fractional integral
		(here indicated as $D_t^{-\nu}$) is not satisfactory because
		\begin{align}
			e^{-\alpha D_t^{-\nu}} = \sum_{r=0}^\infty \frac{(-\alpha)^r D_t^{-r \nu} 1}{r!}
			= \sum_{r=0}^\infty \frac{(-\alpha)^r t^{\nu r}}{r!\Gamma(\nu r+1)}
			= \phi(\nu,1;-\alpha t^\nu),
		\end{align}
		where
		\begin{align}
			\phi(\gamma,\zeta;x) = \sum_{r=0}^\infty \frac{x^r}{r!\Gamma(\gamma r+ \zeta)}
		\end{align}
		is a Wright function (see \citep{kilbas}) and not a Mittag--Leffler function.
		Instead, when $\nu \in (0,1]$, we can operate the right transformation by randomising the time $t$ as follows:
		\begin{align}
			e^{-\alpha t} \longrightarrow \int_0^\infty e^{-\alpha \xi t^\nu} f_{\Xi} (\xi) \mathrm d\xi
			= \mathbb{E}_\Xi e^{-\alpha \Xi t^\nu},
		\end{align}
		or equivalently
		\begin{align}
			\label{acqua}
			t \longrightarrow \log \left[ \mathbb{E}_\Xi e^{-\alpha \Xi t^\nu} \right]^{-1/\alpha},
		\end{align}
		where $\Xi$ is a Wright-distributed random variable, that is with probability density function
		\begin{align}
			f_\Xi(\xi) = \sum_{r=0}^\infty \frac{(-\xi)^r}{r! \Gamma(1-\nu(r+1))}, \qquad \xi > 0, \: \nu \in (0,1].
		\end{align}
		We can therefore write that
		\begin{align}
			\mathbb{E}_\Xi e^{-\alpha \Xi t^\nu} & = \int_0^\infty e^{-\alpha \xi t^\nu} f_\Xi(\xi) \mathrm d\xi \\
			& = \int_0^\infty e^{-\alpha \xi t^\nu} \sum_{r=0}^\infty \frac{(-\xi)^r}{r!\Gamma(1-\nu(r+1))}
			\mathrm d\xi \notag \\
			& = \sum_{r=0}^\infty \frac{(-1)^r}{r! \Gamma(1-\nu r-\nu)}
			\int_0^\infty e^{-\alpha \xi t^\nu} \xi^r \mathrm d \xi \notag \\
			& = \frac{1}{\alpha t^\nu} \sum_{r=0}^\infty \frac{(-1/(\alpha t^\nu))^r}{\Gamma(1-\nu r -\nu)} \notag \\
			& = \frac{1}{\alpha t^\nu} E_{-\nu, 1-\nu}\left( - \frac{1}{\alpha t^\nu} \right) \notag \\
			& = E_{\nu}(-\alpha t^\nu), \notag
		\end{align}
		where
		\begin{align}
			E_{\gamma, \zeta} (y) = \sum_{r=0}^\infty \frac{y^r}{\Gamma(\gamma r + \zeta)}, \qquad
			y \in \mathbb{R},
		\end{align}
		is the generalised Mittag--Leffler function (see \citep{kilbas}), and where, in the last step we
		exploited formula (5.1) of \citep{beg}.
		
		Alternatively, one can define the integral operator
		\begin{align}
			\mathbb{D}_t^{-\nu} f(t) = \Gamma(\lceil \nu \rceil) D_t^{-\nu} f(t), \qquad \nu \in (0,1),
		\end{align}
		which is a simple modification of the Riemann--Liouville fractional integral, and then write that
		\begin{align}
			e^{-\alpha \mathbb{D}_t^{-\nu}} = \sum_{r=0}^\infty \frac{(-\alpha)^r \mathbb{D}_t^{-r\nu} 1}{r!}
			= \sum_{r=0}^\infty \frac{(-\alpha)^r t^{\nu r}}{\Gamma(\nu r+1)} = E_\nu(-\alpha t^\nu).
		\end{align}
		Beside \eqref{acqua}, we also have
		\begin{align}
			t \longrightarrow \mathbb{D}_t^{-\nu}.
		\end{align}
		In practice in our analysis the operator $\mathbb{D}_t^{-\nu}$ replaces the integral operator
		$\Delta_t^{-1}$ in \citet{Dat1}. In the ordinary case we have that $\Delta_t^{-1} 1 = \int_0^t 1 \, \mathrm ds = t$
		while in our case $\mathbb{D}_t^{-\nu} 1 = \int_0^t 1\, \mathrm ds_\nu = t^\nu \Gamma
		(\lceil \nu \rceil)/\Gamma(1+\nu)$, where
		\begin{align}
			\mathrm ds_\nu = \frac{\mathrm ds \, \Gamma (\lceil \nu \rceil)}{\Gamma(\nu) (t-s)^{1-\nu}}.
		\end{align}
		Note that $\mathrm d s_\nu$ can be linked to the interpretation of fractional integration
		furnished by \citet{tarasov} as a measure for fractal media.
		
	\section{Analytic solutions of a class of fractional differential equations by operational methods}  

		\label{solution}
		We here seek to solve analytically a general class of linear partial fractional differential
		equations by operational methods, using the results recalled in the previous section.

		\begin{Theorem}

			\label{ah}
			Consider the following initial value problem (IVP)  
			\begin{align}
				\label{IVP}
				\begin{cases}
					D_{t}^{\nu}f(x,t)=\Theta_x f(x,t), & \nu \in (0,1], \\ 
					f(x,0)=g(x),
				\end{cases}
			\end{align}
			in the half plane $t \ge 0$, with an analytic initial condition $g(x)$ and where $\Theta_x$
			is a generic linear integro-differential operator with constant coefficients acting on x,
			and which satisfies the 
			semi-group property, i.e.\ $\Theta_x \Theta_y = \Theta_{x+y}$.
			The operational solution of equation \eqref{IVP} is given by:
			\begin{align}
				f(x,t)=E_{\nu}(t^{\nu}\Theta_x)g(x)=\sum_{r=0}^{\infty}\frac{t^{\nu r} \Theta^r_x}{\Gamma(r\nu+1)}g(x).
			\end{align}

			\begin{proof}
				Using spectral properties of Caputo fractional derivative, we immediately have that
				\begin{equation}\nonumber
					D_t^{\nu}f(x,t) =  D_t^{\nu}E_{\nu}(t^{\nu}\Theta_x)g(x) =
					\Theta_x E_{\nu}(t^{\nu}\Theta_x)g(x) =  \Theta_x f(x,t).
				\end{equation}
			\end{proof}

		\end{Theorem}

		\begin{Corollary}

			\label{unaparola}
			Consider the following BVP:  
			\begin{align}
				\label{BVP}
				\begin{cases}
					D_{x}^{\nu}f(x,t)= \Theta_t f(x,t), \\
					f(0,t)=g(t),
				\end{cases}
			\end{align}
			in the half plane $x \ge 	0$, with an analytic boundary condition $g(t)$.
			The operational solution of equation \eqref{BVP} is given by:
			\begin{align}
				f(x,t)=E_{\nu}(x^{\nu}\Theta_t)g(t)= \sum_{r=0}^{\infty}\frac{x^{\nu r} \Theta^r_t}{\Gamma(r\nu+1)}g(t).
			\end{align}

		\end{Corollary} 
		
		\begin{Remark}
			Note that Theorem \ref{ah} is still valid for $\nu \in (1,2]$ if the constraint $\partial_t f(x,0)=0$
			holds. Moreover we can state a similar remark for Corollary \ref{unaparola} but with the constraint
			$\partial_x f(0,t)=0$.
		\end{Remark}

		\begin{Ex}[IVP for the time-fractional diffusion equation]

			Consider the following IVP
			\begin{align}
				\begin{cases}
					D_t^{\nu}f(x,t)=\partial_{xx} f(x,t), & \nu \in (0,1],\\
					f(x,0)= x^{\beta}, & \beta \notin \mathbb{Z}^-.
				\end{cases}
			\end{align}
			Considering Theorem \ref{ah}, the analytic solution is given by
			\begin{align}
				\label{computer}
				f(x,t)= \sum_{r=0}^{\infty}\frac{\Gamma(\beta+1)
				t^{\nu r}x^{\beta-2r}}{\Gamma(r\nu+1)\Gamma(\beta+1-2r)}.
			\end{align}
		
			\begin{proof}
				Recalling that $D_t^{\nu}t^{r\nu}=\Gamma(r\nu+1)t^{r\nu-\nu}/\Gamma(r\nu+1-\nu)$,
				we have that
				\begin{align}
					\label{timedif}
					D_t^{\nu} \left\{ \sum_{r=0}^{\infty}\frac{\Gamma(\beta+1)
					t^{\nu r}x^{\beta-2r}}{\Gamma(r\nu+1)\Gamma(\beta+1-2r)}\right\} =
					\sum_{r=1}^{\infty}\frac{\Gamma(\beta+1)t^{-\nu}
					t^{\nu r}x^{\beta-2r}}{\Gamma(r\nu+1-\nu)\Gamma(\beta+1-2r)}.
				\end{align}
				On the other hand
				\begin{align}
					\partial_{xx} \left\{\sum_{r=0}^{\infty}\frac{\Gamma(\beta+1)
					t^{\nu r}x^{\beta-2r}}{\Gamma(r\nu+1)\Gamma(\beta+1-2r)}\right\} =
					\sum_{r=0}^{\infty}\frac{\Gamma(\beta+1)
					t^{\nu r}x^{\beta-2(r+1)}}{\Gamma(r\nu+1)\Gamma(\beta-1-2r)},
				\end{align}
				so with a simple change of index in the sum ($r'= r+1$) we retrieve
				\eqref{timedif} and the proof is complete.

				Note that the solution \eqref{computer} is the fractional generalisation of the
				classical heat polynomials studied by \citet{rose}.

			\end{proof}

		\end{Ex}

		\begin{Ex}[IVP for the time-fractional equation of vibrating plates \citep{mirko}]

			Consider the following initial value problem for the time-fractional equation of vibrating plates
			\begin{align}
				\begin{cases}
					D_{t}^{\nu}f(x,t)=-\partial_{xxxx} f(x,t), & \nu \in (0,2], \\ 
					f(x,0)= \sin x, \\
					\partial_t f(x,0) = 0,
				\end{cases}
			\end{align}
			The analytic solution is given by
			\begin{equation}
				f(x,t)= \sin x \sum_{r=0}^{\infty}\frac{(-1)^{r}t^{\nu r}}{\Gamma(r\nu+1)} = \sin x E_{\nu}(-t^{\nu}).
			\end{equation}
			Note that this is simply the solution by separation of variables.

		\end{Ex}			

		\begin{Ex}[BVP for the space-fractional diffusion equation]

			Consider the following BVP
			\begin{align}
				\begin{cases}
					D_{x}^{\nu}f(x,t)=\partial_{t} f(x,t), & \nu \in (0,2], \: x \ge 0,\\
					f(0,t)= e^{-t},\\
					\partial_x f(0,t) = 0.
				\end{cases}
			\end{align}
			The analytic solution is given by
			\begin{equation}
				f(x,t)= e^{-t} \sum_{r=0}^{\infty}\frac{(-1)^{r}x^{\nu r}}{\Gamma(r\nu+1)}= e^{-t} E_{\nu}(-x^{\nu}).
			\end{equation}
			It is straightforward to realise that this solution can also be derived by separation of variables. 

		\end{Ex}
		
		\begin{Ex}[Fractional Poisson process]
			Consider the differential equation governing the probability generating function $G^\nu(u,t)$,
			$|u|\leq 1$, of a
			fractional Poisson process \citep{beg} $N^\nu(t)$, $t \ge 0$, with rate $\kappa>0$:
			\begin{align}
				\label{mammuth}
				\begin{cases}
					D^\nu_t G^\nu(u,t) = - \kappa (1-u) G^\nu(u,t), \\
					G(u,0) = 1.  
				\end{cases}
			\end{align}
			It may be immediately realised that the solution to \eqref{mammuth} is
			\begin{align}
				G^\nu(u,t) = E_\nu(-\kappa (1-u)t^\nu), \qquad |u| \leq 1.
			\end{align}
			The remarkable thing here is that we can treat the related fractional difference-differential
			equations governing the state probabilities $p_k^\nu(t) = \text{Pr} \{ N^\nu(t) = k \}$, $k \geq 0$,
			with the same operational techniques.
			In particular it is known that $p_k^\nu(t)$ satisfy the equations
			\begin{align}
				\label{ci}
				\begin{cases}
					D^\nu p_k^\nu(t) = - \lambda p_k^\nu(t) + \lambda p_{k-1}^\nu(t), & \nu \in (0,1], \: t \ge 0, \\
					p_k^\nu(0) = \delta_{k,0}.
				\end{cases}
			\end{align}
			Where $\delta_{k,0}$ is the Kronecker delta, i.e.\
			\begin{align}
				\delta_{k,0} =
				\begin{cases}
					1, & k=0, \\
					0, & k\neq 0.
				\end{cases}
			\end{align}
			We can rewrite equation\ \eqref{ci} as
			\begin{align}
				\label{ci2}
				\begin{cases}
					D^\nu p_k^\nu(t) = - \lambda (1-B) p_k^\nu(t), & \nu \in (0,1], \: t \ge 0, \\
					p_k^\nu(0) = \delta_{k,0}.
				\end{cases}
			\end{align}
			where $B$ is the so-called \emph{backward shift operator} and is such that $B(\delta_{k,0})=\delta_{k-1,0}$ and
			$B^r(\delta_{k,0}) = B^{r-1}(B(\delta_{k,0})) = \delta_{k-r,0}$.
			By considering that the eigenfunction associated with the Caputo fractional derivative is the Mittag--Leffler
			function, we obtain that the operational solution to \eqref{ci2} is
			\begin{align}
				\label{qq}
				p_k^\nu(t) & = \sum_{r=0}^\infty \frac{(-\lambda t^\nu)^r}{\Gamma(\nu r+1)} (1-B)^r \delta_{k,0} \\
				& = \sum_{r=0}^\infty \frac{(-\lambda t^\nu)^r}{\Gamma(\nu r+1)} \sum_{h=0}^r
				\binom{r}{h} (-1)^h \delta_{k-h,0} \notag \\
				& = \sum_{r=k}^\infty \frac{(-\lambda t^\nu)^r}{\Gamma(\nu r+1)} \binom{r}{k} (-1)^k. \notag
			\end{align}
			Note that \eqref{qq} coincides with formula (1.4) of \citet{beg}.
			For the specific case $\nu=1$, that is for the classical homogeneous Poisson process, we can retrieve
			the Poisson
			distribution either by specialising \eqref{qq} or by operational methods as follows.
			\begin{align}
				p_k(t) & = \sum_{r=0}^\infty \frac{[-\lambda t(1-B)]^r}{r!} \delta_{k,0} \\
				& = \sum_{r=0}^\infty \frac{(-\lambda t)^r}{r!} \sum_{h=0}^r \binom{r}{h} (-1)^h \delta_{k-h,0} \notag \\
				& = \sum_{r=k}^\infty \frac{(-\lambda t)^r}{r!} \binom{r}{k} (-1)^k \notag \\
				& = \sum_{r=0}^\infty \frac{(-\lambda t)^{r+k}}{k!r!}(-1)^k \notag \\
				& = \frac{(\lambda t)^k}{k!} e^{-\lambda t}. \notag
			\end{align}
						
		\end{Ex}

	\bibliography{bvp} 
	\bibliographystyle{unsrtnat}
	\nocite{*}

\end{document}